\documentclass[12pt]{article}                                               

\input{amssym.def}                                                           
\input{amssym.tex}                                                           
                                                                             
\begin{document}                                                             
\title{On a Cuntz-Krieger functor}

\author{Igor V. Nikolaev
\footnote{Partially supported 
by NSERC.}}


\date{}
 \maketitle


\newtheorem{thm}{Theorem}
\newtheorem{lem}{Lemma}
\newtheorem{dfn}{Definition}
\newtheorem{rmk}{Remark}
\newtheorem{cor}{Corollary}
\newtheorem{prp}{Proposition}
\newtheorem{exm}{Example}
\newtheorem{cnj}{Conjecture}

\begin{abstract}
We construct a covariant functor from the topological torus
bundles to the so-called Cuntz-Krieger algebras; the functor
maps homeomorphic bundles into the stably isomorphic Cuntz-Krieger
algebras. It is shown, that the $K$-theory of the Cuntz-Krieger
algebra encodes torsion of the first homology group of the
bundle. We illustrate the result by examples.

\vspace{7mm}


\vspace{5mm}
{\it AMS (MOS) Subj. Class.: 46L; 57M}
\end{abstract}

\section{Introduction}
{\bf A. The Cuntz-Krieger algebras.}
Recall,  that the Cuntz-Krieger algebra, ${\cal O}_A$, is a  $C^*$-algebra
generated by partial isometries $s_1,\dots, s_n$;  they act on the Hilbert space in such 
a way,  that their support projections $Q_i=s_i^*s_i$ and their  range projections
$P_i=s_is_i^*$ are orthogonal and satisfy the relations
\begin{equation}
Q_i=\sum_{j=i}^na_{ij}P_j,
\end{equation}
for an $n\times n$  matrix $A=(a_{ij})$ consisting of $0$'s 
and $1$'s \cite{CuKr1}.  The notion is extendable to the 
matrices $A$ with the non-negative integer entries {\it ibid., 
Remark 2.16.}  It is known,  that the $C^*$-algebra ${\cal O}_A$ is simple, 
whenever matrix $A$ is irreducible;  the latter means that  certain power
of $A$ is a strictly positive integer matrix.  If ${\cal K}$ is
the $C^*$-algebra of  compact operators on a Hilbert space,
then the Cuntz-Krieger algebras ${\cal O}_A$, ${\cal O}_{A'}$
are said to be {\it stably isomorphic},  if
${\cal O}_A\otimes {\cal K}\cong {\cal O}_{A'}\otimes {\cal K}$,
where $\cong$ is an isomorphism of the $C^*$-algebras.
The $K$-theory of  ${\cal O}_A$  was established in \cite{CuKr1};
it was shown, that $K_0({\cal O}_A)\cong {\Bbb Z}^n / (I-A^t){\Bbb Z}^n$
and $K_1({\cal O}_A)=Ker~(I-A^t)$, where $A^t$ is a transpose of the matrix $A$.
It is easy  to see, that if $det~(I-A^t)\ne 0$, then $K_0({\cal O}_A)$
is a finite abelian group and $K_1({\cal O}_A)=0$;   the two groups are invariants 
of  the  stable isomorphism.

\medskip\noindent
{\bf B. The torus bundles.}
Let $T^n$  be a torus of dimension $n\ge 1$.  Recall, that {\it torus bundle}
is  an  $(n+1)$-dimensional manifold  
\begin{equation}
M_{\alpha}=\{T^n\times [0,1] ~|~ (T^n,0)=(\alpha(T^n),1)\},
\end{equation}
where $\alpha: T^n\to T^n$ is an automorphism of $T^n$. 
It is known, that  bundles $M_{\alpha}$ and $M_{\alpha'}$ are homeomorphic, 
whenever  the automorphisms $\alpha$ and $\alpha'$ are conjugate, 
i.e. $\alpha'=\beta\circ\alpha\circ\beta^{-1}$ for an automorphism 
$\beta: T^n\to T^n$.  Let $H_1(T^n;  {\Bbb Z})\cong {\Bbb Z}^n$ be the first homology of
torus;  consider the group $Aut~(T^n)$ of (homotopy classes of) automorphisms
of $T^n$.   Any $\alpha\in Aut~(T^n)$  induces a linear transformation 
of  $H_1(T^n; {\Bbb Z})$,  given by an invertible $n\times n$  matrix $A$
with the integer entries;  conversely,  each $A\in GL_n({\Bbb Z})$  
defines an automorphism $\alpha: T^n\to T^n$.  In this  matrix 
representation, the conjugate  automorphisms $\alpha$ and  $\alpha'$  define  similar 
matrices  $A,A'\in GL_n({\Bbb Z})$, i.e. such that $A'=BAB^{-1}$ for a matrix 
$B\in GL_n({\Bbb Z})$.  Each class of matrices, similar to  
 a matrix $A\in GL_n({\Bbb Z})$ and such that $tr~(A)\ge 0$ ($tr~(A)\le 0$), 
contains a matrix with only  the non-negative (non-positive) entries.
We always assume, that our bundle $M_{\alpha}$ is given by a non-negative
matrix $A$;  the matrices with  $tr~(A)\le 0$ can be  reduced to this case 
by switching the sign (from negative to positive) in the respective non-positive
representative.

\medskip\noindent
{\bf C. The result.}
Denote by ${\cal M}$ a category of torus bundles (of fixed dimension),  endowed with
homeomorphisms between the bundles; denote by ${\cal A}$ 
a category of the Cuntz-Krieger algebras  ${\cal O}_A$ with  $det~(A)=\pm 1$,
endowed with stable isomorphisms between the algebras.  
Consider a (Cuntz-Krieger) map,  $F: {\cal M}\to {\cal A}$,  which acts by the 
formula $M_{\alpha}\mapsto {\cal O}_A$. The following is true. 
\begin{thm}\label{thm1}
The map $F$ is a covariant functor, such that
$H_1(M_{\alpha}; {\Bbb Z})$ and ${\Bbb Z}\oplus K_0(F(M_{\alpha}))$ 
are isomorphic  abelian groups.
\end{thm}
The article is organized as follows.  Theorem \ref{thm1} is proved in section 2.
There is no formal section on the preliminaries;  all  necessary results
and concepts are introduced in passing. (We encourage the  reader
to consult \cite{CuKr1}--\cite{Wag1} for the details.) 
In section 3,  an application of theorem \ref{thm1} is considered.

\section{Proof of theorem 1}
The idea of proof consists in a reduction of the conjugacy problem for the
automorphisms of $T^n$ to the Cuntz-Krieger theorem on the flow equivalence 
of the subshifts of finite type (to be introduced in the next paragraph).  
There  are no difficult parts in the proof, which is  basically a series of 
observations. Moreover,  theorem \ref{thm1}  follows from the results of P.~M.~Rodrigues and 
J.~S.~Ramos \cite{RoRa1}.  However, our accents are different and the proof is more direct (and shorter)
than in the above cited work.

\medskip
(i) The main reference to the subshifts of finite type (SFT) is \cite{LM}.
A  {\it full} Bernoulli $n$-shift is the set $X_n$ of bi-infinite sequences
$x=\{x_k\}$, where $x_k$ is a symbol taken from a set $S$ of cardinality $n$.  
The set $X_n$ is endowed with the product topology, making $X_n$ a Cantor set.
The shift homeomorphism $\sigma_n:X_n\to X_n$ is given by the formula
$\sigma_n(\dots x_{k-1}x_k x_{k+1}\dots)=(\dots x_{k}x_{k+1}x_{k+2}\dots)$
The homeomorphism defines a (discrete) dynamical system $\{X_n,\sigma_n\}$
given by the iterations of $\sigma_n$.

Let $A$ be an $n\times n$ matrix, whose entries $a_{ij}:= a(i,j)$ are $0$ or $1$.
Consider a subset $X_A$ of $X_n$ consisting of the bi-infinite sequences,
which satisfy the restriction $a(x_k, x_{k+1})=1$ for all $-\infty<k<\infty$.
(It takes a moment to verify that $X_A$ is indeed a subset of $X_n$ and 
$X_A=X_n$, if and only if,  all the entries of $A$ are $1$'s.) By definition,
$\sigma_A=\sigma_n~|~X_A$ and the pair $\{X_A,\sigma_A\}$ is called a {\it SFT}.
A standard edge shift construction described in \cite{LM} allows to extend
the notion of SFT to any matrix $A$ with the non-negative entries.

It is well known that the SFT's $\{X_A,\sigma_A\}$ and $\{X_B,\sigma_B\}$
are topologically conjugate (as the dynamical systems), if and only if, 
the matrices $A$ and $B$ are {\it strong shift equivalent} (SSE), see \cite{LM} for 
the corresponding definition. The SSE of two matrices is a difficult
algorithmic problem, which motivates the consideration of a weaker 
equivalence between the matrices called a {\it shift equivalence} (SE). 
Recall,  that the matrices $A$ and $B$ are said to be shift equivalent
(over ${\Bbb Z}^+$),  when there exist non-negative matrices $R$ and $S$
and a positive integer $k$ (a lag), satisfying the equations 
$AR=RB, BS=SA, A^k=RS$ and $SR=B^k$. Finally, the SFT's    
 $\{X_A,\sigma_A\}$ and $\{X_B,\sigma_B\}$ (and the matrices $A$ and $B$)
are said to be {\it flow equivalent} (FE),  if the suspension flows of the SFT's
act on the topological spaces, which are homeomorphic under a homeomorphism
that respects the orientation of the orbits of the suspension flow. 
We shall use the following implications:
\begin{equation}\label{eq1}
SSE\Rightarrow  SE\Rightarrow  FE.
\end{equation}
(The first implication is rather classical, while for the second
we refer the reader to \cite{LM}, p.456.)

We further restrict to the SFT's given by the matrices with determinant
$\pm 1$. In view of Corollary 2.13 of \cite{Wag1},  the matrices $A$ and
$B$ with $det~(A)=\pm 1$ and $det~(B)=\pm 1$ are SE (over ${\Bbb Z}^+$),
if and only if,  matrices $A$ and $B$ are similar  in $GL_n({\Bbb Z})$.

Let now $\alpha$ and $\alpha'$ be a pair of conjugate automorphisms of $T^n$.
Since the corresponding matrices $A$ and $A'$ are similar in $GL_n({\Bbb Z})$,
one concludes that the SFT's $\{X_{A},\sigma_{A}\}$ and 
$\{X_{A'},\sigma_{A'}\}$ are SE. In particular, 
the SFT's $\{X_{A},\sigma_{A}\}$ and 
$\{X_{A'},\sigma_{A'}\}$ are FE.

One can now apply the known result due to Cuntz and Krieger; it says, 
that the $C^*$-algebra ${\cal O}_A\otimes {\cal K}$ is an invariant
of the flow equivalence of the irreducible SFT's, see p. 252 of \cite{CuKr1}
and its proof in Sect. 4 of the same work.  Thus, the map $F$ sends the 
conjugate automorphisms of $T^n$ into the stably
isomorphic Cuntz-Krieger algebras, i.e. $F$ is a functor.

Let us show that $F$ is a covariant functor. Consider the following commutative
diagram:

\begin{picture}(300,110)(-80,-5)
\put(20,70){\vector(0,-1){35}}
\put(130,70){\vector(0,-1){35}}
\put(45,23){\vector(1,0){53}}
\put(45,83){\vector(1,0){53}}
\put(17,20){${\cal O}_{A}$}
\put(5,55){$F$}
\put(140,55){$F$}
\put(118,20){${\cal O}_{BAB^{-1}},$}
\put(17,80){$A$}
\put(117,80){$A'=BAB^{-1}$}
\put(55,30){\sf stable}
\put(50,12){\sf isomorphism}
\put(54,90){\sf similarity}
\end{picture}

\noindent
where $A,B\in GL_n({\Bbb Z})$ and ${\cal O}_A, {\cal O}_{BAB^{-1}}\in {\cal A}$.  
Let $g_1,g_2$ be the arrows (similarity
of matrices) in the upper category and $F(g_1), F(g_2)$ the corresponding
arrows (stable  isomorphisms) in the lower category. In view of the
diagram, we have the following identities:
\begin{eqnarray}
F(g_1g_2) &=& {\cal O}_{B_2B_1AB_1^{-1}B_2^{-1}}= {\cal O}_{B_2(B_1AB_1^{-1})B_2^{-1}}\nonumber\\
          &=& {\cal O}_{B_2A'B_2^{-1}}=F(g_1)F(g_2),
\end{eqnarray}
where $F(g_1)({\cal O}_A)={\cal O}_{A'}$ and $F(g_2)({\cal O}_{A'})={\cal O}_{A''}$.
Thus, $F$ does not reverse the arrows and is,  therefore,  a covariant functor.  
The first statement of theorem \ref{thm1} is proved.

\bigskip
(ii) Let $M_{\alpha}$ be a torus bundle with a monodromy,  given by the matrix $A\in GL_n({\Bbb Z})$.
It can be calculated, e.g. using the Leray spectral sequence for the fiber bundles, 
that $H_1(M_{\alpha}; {\Bbb Z})\cong {\Bbb Z}\oplus  {\Bbb Z}^n / (A-I){\Bbb Z}^n$. 
Comparing this  calculation with the $K$-theory of the Cuntz-Krieger algebra, one concludes
that $H_1(M_{\alpha}; {\Bbb Z})\cong {\Bbb Z}\oplus K_0({\cal O}_{A})$, where ${\cal O}_{A}=F(M_{\alpha})$. 
The second statement   of theorem \ref{thm1} follows.
$\square$

\section{Examples}
Consider  the following (three-dimensional) torus bundles $M_{\alpha_i} ~(i=1,2,3)$
\footnote{It is interesting,  that by  Thurston's Geometrization Theorem,  the bundle $M_{\alpha_1^n}$ is  a 
nilmanifold for any  $n$, while bundles $M_{\alpha_2}$ and  $M_{\alpha_3}$ are the 
solvmanifolds \cite{Thu1}.}
:
\begin{eqnarray}
A_1^n &=& \left(\matrix{1 & n\cr 0 & 1}\right), \quad n\in {\Bbb Z}, 
\qquad K_0({\cal O}_{A_1^n})\cong {\Bbb Z}\oplus {\Bbb Z}_n,\nonumber\\
A_2 &=& \left(\matrix{5 & 2\cr 2 & 1}\right),
\qquad K_0({\cal O}_{A_2})\cong {\Bbb Z}_2\oplus {\Bbb Z}_2,\nonumber\\
A_3 &=& \left(\matrix{5 & 1\cr 4 & 1}\right),
\qquad K_0({\cal O}_{A_3})\cong {\Bbb Z}_4. \nonumber
\end{eqnarray}
(Here the group  $K_0({\cal O}_{A_i})$ was calculated using a reduction of 
the matrix  to its Smith normal form \cite{LM}.) 
Notice, that the Cuntz-Krieger invariant  $K_0({\cal O}_{A_1^n})\cong {\Bbb Z}\oplus {\Bbb Z}_n$
is a complete topological invariant  of the family of bundles  $M_{\alpha_1^n}$; 
thus, such an invariant solves a classification problem for these bundles.
For the bundles $M_{\alpha_2}$ and $M_{\alpha_3}$ the Alexander polynomial:
\begin{equation}
\Delta_{A_2}(t)=\Delta_{A_3}(t)=t^2-6t+1.
\end{equation}
Thus,  the Alexander polynomial alone   cannot  distinguish between the bundles
$M_{\alpha_2}$ and $M_{\alpha_3}$;  however, since   
$K_0({\cal O}_{A_2})\not\cong  K_0({\cal O}_{A_3})$,  theorem \ref{thm1} says 
that the bundles $M_{\alpha_2}$ and $M_{\alpha_3}$ are topologically distinct.

\bigskip\noindent
{\sf Acknowledgments.} 
The paper is dedicated to Prof. J. ~Cuntz on occasion of his sixtieth birthday.  
I am grateful to Prof. W. ~Krieger for a copy of \cite{RoRa1}
and exciting discussions.



\vskip1cm

\textsc{The Fields Institute for Mathematical Sciences, Toronto, ON, Canada,  
E-mail:} {\sf igor.v.nikolaev@gmail.com}

\smallskip
{\it Current address: 101-315 Holmwood Ave., Ottawa, ON, Canada, K1S 2R2}


\begin{thebibliography}{100}
\bibitem{CuKr1}
J. ~Cuntz and W. ~Krieger, A class of $C^*$-algebras and topological Markov
chains, Invent. Math. 56 (1980), 251-268. 


\bibitem{LM}
D.~Lind and B.~Marcus, An Introduction to Symbolic Dynamics and Coding,
Cambridge Univ. Press, 1995. 





\bibitem{RoRa1}
P.~M.~Rodrigues and J.~S.~Ramos,  Bowen-Franks groups as conjugacy invariants
for $T^n$-automorphisms,  Aequationes Math. 69 (2005), 231-249.  


\bibitem{Thu1}
W.~P.~Thurston, Three dimensional manifolds, Kleinian groups and
hyperbolic geometry, Bull. Amer. Math. Soc. 6 (1982), 357-381.




\bibitem{Wag1}
J.~B.~Wagoner, Strong shift equivalence theory and the shift equivalence
problem, Bull. Amer. Math. Soc. 36 (1999), 271-296. 






\end{thebibliography}
\end{document}